# EXOTIC SERIES WITH BERNOULLI, HARMONIC, CATALAN, AND STIRLING NUMBERS

Khristo N. Boyadzhiev

Department of Mathematics
Ohio Northern University
Ada, OH 45810, USA
k-boyadzhiev@onu.edu

**Abstract**

In this paper, we present a formula for generating various "exotic" series in the spirit of Ovidiu Furdui and Alina Sintamarian [5]. Our new series (evaluated in closed form) involve Bernoulli, harmonic, and Catalan numbers. Also Stirling numbers of the second kind, other special numbers, and exponential polynomials. The results include series identities with Laguerre polynomials and derangement polynomials.

## 1. Introduction,

In a recent paper Ovidiu Furdui and Alina Sintamarian [5] evaluated several interesting exotic series. For example, they proved that

$$(1) \qquad \sum_{n=1}^{\infty}\left(e^x - 1 - \frac{1}{1!} - \frac{1}{2!} - \ldots - \frac{1}{n!}\right)x^n = \frac{e^x - ex}{x-1} + 1$$

for every $x \neq 1$ with a limit case for $x = 1$

$$\sum_{n=1}^{\infty}\left(e^x - 1 - \frac{1}{1!} - \frac{1}{2!} - \ldots - \frac{1}{n!}\right) = 1.$$

They also proved the evalation

(2) $$\sum_{n=k}^{\infty}\binom{n}{k}\left(e - 1 - \frac{1}{1!} - \ldots - \frac{1}{n!}\right)x^n = \frac{x^k e}{(1-x)^{k+1}}\left\{1 - e^{-(1-x)}\sum_{j=0}^{k}\frac{(1-x)^j}{j!}\right\}$$

for $x \neq 1$ with the limit case for $x = 1$

$$\sum_{n=k}^{\infty}\binom{n}{k}\left(e - 1 - \frac{1}{1!} - \ldots - \frac{1}{n!}\right) = \frac{e}{(k+1)!}.$$

The series (2) appeared as Monthly Problem 12012 (Amer. Math. Monthly, vol. 124, December 2017) proposed by the same authors.

Furdui in another publication [6] evaluated the series

(3) $$\sum_{n=1}^{\infty} n^p \left(e^y - 1 - \frac{y}{1!} - \frac{y^2}{2!} - \ldots - \frac{y^n}{n!}\right) \quad (p \geq 1)$$

which was discussed later by the present author in [2].

In this paper we develop a unified approach for the construction and evaluation of such series. Our general result includes the cases (1), (2), and (3) and produces further interesting exotics series. The main theorem is given in the next section and in section 3 we present the applications. For illustration, three exotic series proved in Section 3 are

$$\sum_{n=0}^{\infty} B_n \left(e^y - 1 - \frac{y}{1!} - \ldots - \frac{y^n}{n!}\right)(-1)^n = e^y\left(y\ln(1-e^{-y}) - \operatorname{Li}_2(e^{-y}) + \frac{\pi^2}{6}\right)$$

$$\sum_{n=0}^{\infty} H_n \left(e^y - 1 - \frac{y}{1!} - \ldots - \frac{y^n}{n!}\right) = e^y(y\operatorname{Ein}(y) - y + 1) - 1$$

$$\sum_{n=0}^{\infty} \varphi_n(\lambda) \left(e^y - 1 - \frac{y}{1!} - \ldots - \frac{y^n}{n!}\right) = \frac{e^y}{\lambda}(1 - e^{\lambda(e^{-y}-1)})$$

where $B_n$ are the Bernoulli numbers, $H_n$ are the harmonic numbers, $\text{Ein}(z)$ is the exponential integral function (see below equation (24)), and $\varphi_n(\lambda)$ are the exponential polynomials. At the end we also prove series identities involving derangement polynomials and Laguerre polynomials (examples 8 and 10).

## 2. Generating exotic series

Suppose we have a function $F(z)$ analytic in a neighborhood of the origin and written in the form

(4) $$F(z) = \sum_{n=0}^{\infty} a_n \frac{z^n}{n!}.$$

We construct a function of two variables using the coefficients of $F(z)$

(5) $$u(x,y) = \sum_{n=0}^{\infty} a_n \left( e^y - 1 - \frac{y}{1!} - \ldots - \frac{y^n}{n!} \right) x^n = a_0(e^y - 1) + a_1 \left( e^y - 1 - \frac{y}{1!} \right) x + \ldots \,.$$

To better understand the convergence of the above series we show a rough estimate for its terms. Let $a > 0$ be arbitrary. For $|y| < a$ we have the following estimate from Taylor's formula with reminder in the form of Lagrange

$$\left| e^y - 1 - \frac{y}{1!} - \ldots - \frac{y^n}{n!} \right| \leq \frac{e^a |y|^{n+1}}{(n+1)!}$$

(here $a$ does not depend on $n$). From this

(6) $$\left| a_n \left( e^y - 1 - \frac{y}{1!} - \ldots - \frac{y^n}{n!} \right) x^n \right| \leq \frac{e^a |a_n||y|^{n+1}|x|^n}{(n+1)!} = \frac{|a_n||x|^n}{n!} \frac{|y|^{n+1} e^a}{n+1}.$$

We will keep this estimate in mind when we consider the applications.

**Theorem.** Let $F(z)$ and $u(x,y)$ be as in (4) and (5). Then for all appropriate value of $x, y$ we have the integral representation

(7) $$u(x, y) = e^y \int_0^y e^{-t} F(xt) \, dt.$$

*Proof.* We compute the partial derivative of $u(x, y)$ with respect to $y$

$$u_y(x, y) = \sum_{n=0}^{\infty} a_n \left( e^y - 1 - \frac{y}{1!} - \ldots - \frac{y^{n-1}}{(n-1)!} \right) x^n.$$

Then clearly

$$u_y - u = \sum_{n=0}^{\infty} \frac{a_n}{n!} y^n x^n = F(xy).$$

For $x$ fixed this is a linear differential equation (with respect to the variable $y$) with integrating factor $e^{-y}$. That is,

$$\frac{d}{dy}(ue^{-y}) = e^{-y} F(xy).$$

From here

$$ue^{-y} = \int_0^y e^{-t} F(xt) \, dt + C(x)$$

where $C(x)$ is the constant of integration with respect to $y$. With $y = 0$ we find $C(x) = u(x, 0) = 0$. This way

$$u(x, y) = e^y \int_0^y e^{-t} F(xt) \, dt$$

and the theorem is proved.

The theorem can be proved also by using Taylor's formula with integral remainder.

**Example 1**. Taking the function

$$F(z) = e^z = \sum_{n=0}^{\infty} \frac{z^n}{n!}, \quad a_n = 1$$

we find

$$u(x, y) = \sum_{n=0}^{\infty}\left(e^y - 1 - \frac{y}{1!} - \ldots - \frac{y^n}{n!}\right)x^n = e^y \int_0^y e^{-t} e^{xt}\, dt.$$

Then for $x = 1$ we have

(8) $$\sum_{n=0}^{\infty}\left(e^y - 1 - \frac{y}{1!} - \ldots - \frac{y^n}{n!}\right) = y e^y$$

and for $x \neq 1$

(9) $$\sum_{n=0}^{\infty}\left(e^y - 1 - \frac{y}{1!} - \ldots - \frac{y^n}{n!}\right)x^n = \frac{e^{xy} - e^y}{x - 1}$$

which confirms (1). (Note that in (1) the summation starts from $n = 1$.)

**Example 2.** Now let $p \geq 0$ be an integer. Taking the exponential generating function for the binomial coefficient

$$F(z) = \frac{z^p e^z}{p!} = \sum_{n=0}^{\infty}\binom{n}{p}\frac{z^n}{n!}$$

we find from (7)

$$\sum_{n=0}^{\infty}\binom{n}{p}\left(e^y - 1 - \frac{y}{1!} - \ldots - \frac{y^n}{n!}\right)x^n = e^y \int_0^y e^{-t}\left\{\frac{(xt^p)e^{xt}}{p!}\right\}dt = \frac{e^y x^p}{p!}\int_0^y t^p e^{(x-1)t}\, dt.$$

For $x = 1$

(10) $$\sum_{n=0}^{\infty}\binom{n}{p}\left(e^y - 1 - \frac{y}{1!} - \ldots - \frac{y^n}{n!}\right) = \frac{e^y y^{p+1}}{(p+1)!}$$

and for $x \neq 1$ integration by parts gives

(11) $$\sum_{n=0}^{\infty}\binom{n}{p}\left(e^y - 1 - \frac{y}{1!} - \ldots - \frac{y^n}{n!}\right)x^n = \frac{x^p e^y}{(1-x)^{p+1}}\left\{1 - e^{-(1-x)y}\sum_{j=0}^{p}\frac{(1-x)^j y^j}{j!}\right\}$$

which confirms Furdui and Sintamarian's result (2) by setting $y = 1$.

**Example 3**. For the series in (3) we consider the function

$$F(z) = \sum_{n=0}^{\infty} n^p \frac{x^n}{n!} = \varphi_p(x)e^x$$

where $p \geq 0$ is an integer, $\varphi_p(x)$ is the exponential polynomial or order $p$

(12) $$\varphi_p(x) = \sum_{k=0}^{p} S(p,k)x^k$$

and $S(p,k)$ are the Stirling numbers of the second kind [1, 3, 7]. Applying the theorem we find the representation (with the notational agreement $0^0 = 1$)

(13) $$\sum_{n=0}^{\infty} n^p \left( e^y - 1 - \frac{y}{1!} - \frac{y^2}{2!} - \ldots - \frac{y^n}{n!} \right) x^n = e^y \int_0^y e^{-t} e^{xt} \varphi_p(xt) dt.$$

When $x = 1$ this becomes

(14) $$\sum_{n=0}^{\infty} n^p \left( e^y - 1 - \frac{y}{1!} - \frac{y^2}{2!} - \ldots - \frac{y^n}{n!} \right) = e^y \int_0^y \varphi_p(t) dt$$

$$= e^y \sum_{k=0}^{p} S(p,k) \frac{y^{k+1}}{k+1}$$

which is Furdui's result [5]. For $x \neq 0, 1$ we have

(15) $$\sum_{n=0}^{\infty} n^p \left( e^y - 1 - \frac{y}{1!} - \frac{y^2}{2!} - \ldots - \frac{y^n}{n!} \right) x^n$$

$$= \sum_{k=0}^{p} S(p,k) \frac{k! x^k}{(1-x)^{k+1}} \left( e^y - e^{xy} \sum_{j=0}^{k} \frac{y^j (1-x)^j}{j!} \right)$$

(see also [2]).

# 3. Further exotic series

In this section we use our theorem to evaluate various exotic series with special numbers.

**Example 4.** Consider the generating function for the Bernoulli numbers

$$F(z) = \frac{z}{e^z - 1} = \sum_{n=0}^{\infty} B_n \frac{z^n}{n!} \quad (|z| < 2\pi).$$

The theorem implies the representation

(16) $$\sum_{n=0}^{\infty} B_n \left( e^y - 1 - \frac{y}{1!} - \ldots - \frac{y^n}{n!} \right) x^n = e^y \int_0^y \frac{e^{-t} xt}{e^{xt} - 1} dt = xe^y \int_0^y \frac{te^{-t}}{e^{xt} - 1} dt.$$

For $x, y > 0$ the integral can be evaluated in terms of series

$$\int_0^y \frac{te^{-t}}{e^{xt} - 1} dt = \int_0^y \frac{te^{-t} e^{-xt}}{1 - e^{-xt}} dt = \int_0^y t \left\{ \sum_{n=1}^{\infty} e^{-(nx+1)t} \right\} dt = \sum_{n=1}^{\infty} \int_0^y te^{-(nx+1)t} dt$$

$$= -ye^{-y} \sum_{n=1}^{\infty} \frac{e^{-nxy}}{nx+1} - e^{-y} \sum_{n=1}^{\infty} \frac{e^{-nxy}}{(nx+1)^2} + \sum_{n=1}^{\infty} \frac{1}{(nx+1)^2}.$$

These series can be expressed through the Lerch transcendent $\Phi(z, s, a)$ (see [4]) and the dilogarithm $\text{Li}_2(z)$

$$\Phi(z, s, a) = \sum_{n=0}^{\infty} \frac{z^n}{(n+a)^s}; \quad \text{Li}_2(z) = \sum_{n=1}^{\infty} \frac{z^n}{n^2}.$$

Namely, we have

(17) $$\sum_{n=0}^{\infty} B_n \left( e^y - 1 - \frac{y}{1!} - \ldots - \frac{y^n}{n!} \right) x^n = \frac{e^y}{x} \Phi(1, 2, x^{-1}) - y \Phi(e^{-xy}, 1, x^{-1}) - \frac{1}{x} \Phi(e^{-xy}, 2, x^{-1})$$

$$+ x(y + 1 - e^y).$$

For $x = 1$ this representation takes the form

(18) $$\sum_{n=0}^{\infty} B_n \left( e^y - 1 - \frac{y}{1!} - \ldots - \frac{y^n}{n!} \right) = e^y \left( y \ln(1 - e^{-y}) - \text{Li}_2(e^{-y}) + \frac{\pi^2}{6} - 1 \right) + 1 + y.$$

For $x = -1$ the integral in (16) becomes simpler

$$-\int_0^y \frac{te^{-t}}{e^{-t}-1} dt = \int_0^y t\, d\ln(1-e^{-t}) = y\ln(1-e^{-y}) - \int_0^y \ln(1-e^{-t})\, dt$$

$$= y\ln(1-e^{-y}) + \frac{\pi^2}{6} - \text{Li}_2(e^{-y})$$

and we come to the remarkable evaluation.

(19) $$\sum_{n=0}^{\infty} B_n\left(e^y - 1 - \frac{y}{1!} - \dots - \frac{y^n}{n!}\right)(-1)^n = e^y\left(y\ln(1-e^{-y}) - \text{Li}_2(e^{-y}) + \frac{\pi^2}{6}\right).$$

**Example 5.** In this example we use the exponential generating function for the Stirling numbers $S(n,k)$ of the second kind [3, 7]

$$F(z) = \frac{1}{k!}(e^z - 1)^k = \sum_{n=0}^{\infty} S(n,k)\frac{z^n}{n!}$$

where $k \geq 0$ is an integer (the summation actually starts from $n = k$, as $S(n,k) = 0$ for $n < k$). From the theorem

(20) $$\sum_{n=0}^{\infty} S(n,k)\left(e^y - 1 - \frac{y}{1!} - \dots - \frac{y^n}{n!}\right)x^n = \frac{e^y}{k!}\int_0^y e^{-t}(e^{xt} - 1)^k\, dt.$$

The integral can be evaluated in terms of binomial expressions.

$$\int_0^y e^{-t}(e^{xt}-1)^k\, dt = \int_0^y \left\{\sum_{j=0}^k \binom{k}{j}(-1)^j e^{(jx-1)t}\right\} dt = \sum_{j=0}^k \binom{k}{j}(-1)^j \int_0^y e^{(jx-1)t}\, dt$$

$$= \sum_{j=0}^k \binom{k}{j}(-1)^j \int_0^y e^{(jx-1)t}\, dt = \sum_{j=0}^k \binom{k}{j}(-1)^j \frac{e^{(jx-1)y}-1}{jx-1}.$$

Thus we have the closed form evaluation

(21) $$\sum_{n=0}^{\infty} S(n,k)\left(e^y - 1 - \frac{y}{1!} - \dots - \frac{y^n}{n!}\right)x^n = \frac{1}{k!}\sum_{j=0}^k \binom{k}{j}(-1)^j \frac{e^{jxy} - e^y}{jx-1}.$$

**Example 6.** Now we construct an exotic series containing the exponential polynomials $\varphi_n(x)$ defined in (12). Their generating function is given by

$$F(z) = e^{x(e^z - 1)} = \sum_{n=0}^{\infty} \varphi_n(x) \frac{z^n}{n!}$$

(see [1]). The theorem implies the representation

(22) $$\sum_{n=0}^{\infty} \varphi_n(\lambda) \left( e^y - 1 - \frac{y}{1!} - \ldots - \frac{y^n}{n!} \right) x^n = e^y \int_0^y e^{-t} e^{\lambda(e^{xt} - 1)} \, dt = e^{y - \lambda} \int_0^y e^{-t} e^{\lambda e^{xt}} \, dt.$$

We can easily evaluate this integral when $x = -1$

$$\int_0^y e^{-t} e^{\lambda e^{-t}} \, dt = -\frac{1}{\lambda} \int_0^y e^{\lambda e^{-t}} \, d\lambda e^{-t} = -\frac{1}{\lambda} e^{\lambda e^{-t}} \Big|_0^y = -\frac{1}{\lambda} (e^{\lambda e^{-y}} - e^{\lambda}).$$

From this

$$e^{y - \lambda} \int_0^y e^{-t} e^{\lambda e^{xt}} \, dt = -\frac{e^y}{\lambda} (e^{\lambda(e^{-y} - 1)} - 1)$$

and we come to the elegant formula

(23) $$\sum_{n=0}^{\infty} \varphi_n(\lambda) \left( e^y - 1 - \frac{y}{1!} - \ldots - \frac{y^n}{n!} \right) = \frac{e^y}{\lambda} (1 - e^{\lambda(e^{-y} - 1)}).$$

Unexpectedly, we can spot the generating function for the exponential polynomials in the expression on the right hand side. This gives the identity

$$\sum_{n=0}^{\infty} \varphi_n(\lambda) \left( e^y - 1 - \frac{y}{1!} - \ldots - \frac{y^n}{n!} \right) = \frac{e^y}{\lambda} \left( 1 - \sum_{n=0}^{\infty} \varphi_n(\lambda) \frac{(-1)^n y^n}{n!} \right).$$

**Example 7.** I this application we present an exotic series with harmonic number. The harmonic numbers are defined by

$$H_n = 1 + \frac{1}{2} + \ldots + \frac{1}{n} \quad (n \geq 1); \; H_0 = 0$$

with exponential generating function

$$F(z) = \sum_{n=0}^{\infty} H_n \frac{z^n}{n!} = e^z \operatorname{Ein}(z).$$

Here $\mathrm{Ein}(z)$ is the exponential integral

(24) $$\mathrm{Ein}(z) = \int_0^z \frac{1-e^{-u}}{u} du = \sum_{n=1}^{\infty} \frac{(-1)^{n-1} z^n}{n!n}$$

(as the series representation shows, $\mathrm{Ein}(z)$ is an entire function).

From the theorem we obtain

(25) $$\sum_{n=0}^{\infty} H_n\left(e^y - 1 - \frac{y}{1!} - \ldots - \frac{y^n}{n!}\right) x^n = e^y \int_0^y e^{-t} e^{xt} \mathrm{Ein}(xt)\, dt.$$

The integral can easily be evaluated in explicit form when $x = 1$

$$\int_0^y e^{-t} e^t \mathrm{Ein}(t)\, dt = \int_0^y \left\{ \int_0^t \frac{1-e^{-u}}{u} du \right\} dt = \int_0^y \left\{ \int_u^y dt \right\} \frac{1-e^{-u}}{u} du$$

$$= \int_0^y (y-u) \frac{1-e^{-u}}{u} du = y\,\mathrm{Ein}(y) - \int_0^y (1-e^{-u})\, du = y\,\mathrm{Ein}(y) - y - e^{-y} + 1.$$

Using the series representation of $\mathrm{Ein}(t)$ we find also

$$\int_0^y \mathrm{Ein}(t)\, dt = \sum_{n=1}^{\infty} \frac{(-1)^{n-1} y^{n+1}}{n!n(n+1)}.$$

Finally, for all $y$ we have the beautiful equation

(26) $$\sum_{n=0}^{\infty} H_n\left(e^y - 1 - \frac{y}{1!} - \ldots - \frac{y^n}{n!}\right) = e^y(y\,\mathrm{Ein}(y) - y + 1) - 1$$

$$= \sum_{n=1}^{\infty} \frac{(-1)^{n-1} y^{n+1}}{n!n(n+1)}.$$

**Example 8.** The representation (7) leads to some interesting identities involving derangement numbers and polynomials. We come to these identities by using the simple exponential generating function for the numbers $n!$

$$F(z) = \sum_{n=0}^{\infty} n! \frac{z^n}{n!} = \sum_{n=0}^{\infty} z^n = \frac{1}{1-z} \quad (|z|<1).$$

Our theorem gives the representation

(27) $$\sum_{n=0}^{\infty} n! \left( e^y - 1 - \frac{y}{1!} - \ldots - \frac{y^n}{n!} \right) x^n = e^y \int_0^y \frac{e^{-t}}{1-xt} dt.$$

This result can be related to the derangement numbers

$$D_n = \sum_{k=0}^{n} \binom{n}{k} (-1)^{n-k} k! = (-1)^n \sum_{k=0}^{n} \binom{n}{k} (-1)^k k!$$

$$= n! \sum_{k=0}^{n} \frac{(-1)^k}{k!} = n! \left( 1 - \frac{1}{1!} + \frac{1}{2!} + \ldots + \frac{(-1)^n}{n!} \right)$$

which are popular in combinatorics [3, p. 180], [7, pp. 194-196], [8, 9]. We have from their definition

$$e^{-1} n! - D_n = n! \left( e^{-1} - 1 + \frac{1}{1!} - \frac{1}{2!} + \ldots + \frac{(-1)^n}{n!} \right)$$

and (27) with $y = -1$ gives the representation

(28) $$\sum_{n=0}^{\infty} (e^{-1} n! - D_n) x^n = e^{-1} \int_0^{-1} \frac{e^{-t}}{1-xt} dt.$$

Now consider also the derangement polynomials [9]

(29) $$d_n(x) = (-1)^n \sum_{k=0}^{n} \binom{n}{k} (-1)^k k! x^k = n! \sum_{j=0}^{n} \frac{(-1)^j}{j!} x^{n-j}$$

where $d_n(1) = D_n$. The exponential generating function for these polynomials can be computed easily

$$\sum_{n=0}^{\infty} d_n(x) \frac{z^n}{n!} = \sum_{n=0}^{\infty} \frac{z^n}{n!} \left\{ (-1)^n \sum_{k=0}^{n} \binom{n}{k} k! (-1)^k x^k \right\}$$

$$= \sum_{k=0}^{\infty}(-1)^k k! x^k \left\{ \sum_{n=k}^{\infty} \binom{n}{k} \frac{(-1)^n z^n}{n!} \right\} = \sum_{k=0}^{\infty}(-x)^k (-t)^k \left\{ \sum_{n=k}^{\infty} \frac{(-z)^{n-k}}{(n-k)!} \right\}$$

$$= \sum_{k=0}^{\infty}(xz)^k \left\{ \sum_{m=0}^{\infty} \frac{(-z)^m}{m!} \right\} = \frac{e^{-z}}{1-xz}.$$

That is,

(30) $$\frac{e^{-z}}{1-xz} = \sum_{n=0}^{\infty} d_n(x) \frac{z^n}{n!}$$

This function appears in (27) and (28). We have by integrating in (30) the representation

$$\int_0^y e^{-t} \frac{1}{1-xt} dt = \sum_{n=0}^{\infty} d_n(x) \frac{y^{n+1}}{(n+1)!}$$

which in view of (27) gives the curious series identity

(31) $$\sum_{n=0}^{\infty} n! \left( e^y - 1 - \frac{y}{1!} - \dots - \frac{y^n}{n!} \right) x^n = e^y \sum_{n=0}^{\infty} d_n(x) \frac{y^{n+1}}{(n+1)!}.$$

For $x=1, \ y=-1$ this identity becomes

(32) $$\sum_{n=0}^{\infty} n! \left( e^{-1} - 1 + \frac{1}{1!} - \dots - \frac{(-1)^n}{n!} \right) = e^{-1} \sum_{n=0}^{\infty} D_n \frac{(-1)^{n+1}}{(n+1)!}$$

$$= \sum_{n=0}^{\infty}(e^{-1} n! - D_n) = e^{-1} \int_0^{-1} \frac{e^{-t}}{1-t} dt.$$

**Example 9.** In this example we construct exotic series with central binomial coefficients $\binom{2n}{n}$ and Catalan numbers $C_n = \binom{2n}{n} \frac{1}{n+1}$. The Catalan numbers especially are very popular in combinatorics ([3, p. 53] and [7, pp. 203 and 358]).

The exponential generating functions for these numbers are

$$\sum_{n=0}^{\infty}\binom{2n}{n}\frac{z^{n}}{n!}=e^{2z}\sum_{n=0}^{\infty}\frac{z^{2n}}{(n!)^{2}}=e^{2z}I_{0}(2z)$$

$$\sum_{n=0}^{\infty}C_{n}\frac{z^{n}}{n!}=\sum_{n=0}^{\infty}\binom{2n}{n}\frac{z^{n}}{(n+1)!}=e^{2z}\left(I_{0}(2z)-I_{1}(2z)\right)$$

where $I_0(x)$ and $I_1(x)=I_0'(x)$ are the modified Bessel functions of the first kind [11, pp. 77-84]

(33) $$I_{0}(z)=\sum_{n=0}^{\infty}\frac{z^{2n}}{4^{n}(n!)^{2}},\ I_{1}(z)=\frac{z}{2}\sum_{n=0}^{\infty}\frac{z^{2n}}{4^{n}(n!)^{2}(n+1)}.$$

The theorem implies

(34) $$\sum_{n=0}^{\infty}\binom{2n}{n}\left(e^{y}-1-\frac{y}{1!}-\ldots-\frac{y^{n}}{n!}\right)x^{n}=e^{y}\int_{0}^{y}e^{(2x-1)t}I_{0}(2xt)\,dt$$

(35) $$\sum_{n=0}^{\infty}C_{n}\left(e^{y}-1-\frac{y}{1!}-\ldots-\frac{y^{n}}{n!}\right)x^{n}=e^{y}\int_{0}^{y}e^{(2x-1)t}(I_{0}(2xt)-I_{0}'(2xt))\,dt\ .$$

For $x=\dfrac{1}{2}$ the integrals become simpler and can be evaluated explicitly. Thus

$$\int_{0}^{y}I_{0}(t)\,dt=yI_{0}(y)+\frac{\pi y}{2}[I_{0}(y)\mathrm{L}_{1}(y)-I_{1}(y)\mathrm{L}_{0}(y)]$$

where $\mathrm{L}_0(y)$, $\mathrm{L}_1(y)$ are the modified Struve functions [10, entry 1.11.1(4)]. We come to the identities

(36) $$\sum_{n=0}^{\infty}\binom{2n}{n}\frac{1}{2^{n}}\left(e^{y}-1-\frac{y}{1!}-\ldots-\frac{y^{n}}{n!}\right)=e^{y}\sum_{n=0}^{\infty}\frac{y^{2n+1}}{4^{n}(n!)^{2}(2n+1)}$$

$$=e^{y}\left(yI_{0}(y)+\frac{\pi y}{2}[I_{0}(y)\mathrm{L}_{1}(y)-I_{1}(y)\mathrm{L}_{0}(y)]\right)$$

(37) $$\sum_{n=0}^{\infty}\frac{C_{n}}{2^{n}}\left(e^{y}-1-\frac{y}{1!}-\ldots-\frac{y^{n}}{n!}\right)=e^{y}\left(1-I_{0}(y)+\sum_{n=0}^{\infty}\frac{y^{2n+1}}{4^{n}(n!)^{2}(2n+1)}\right)$$

$$= e^y \left(1+(y-1)I_0(y)+\frac{\pi y}{2}[I_0(y)L_1(y)-I_1(y)L_0(y)]\right).$$

**Example 10**. Equation (9) can be viewed as the ordinary generating function for the functions

$$e^y - 1 - \frac{y}{1!} - \ldots - \frac{y^n}{n!}$$

We will show now that the exponential generating function for these expressions is very close to the exponential generating function for the Laguerre polynomials $L_n(x)$. Namely, the following series identity holds

(38) $$\sum_{n=0}^{\infty}\left(e^y - 1 - \frac{y}{1!} - \ldots - \frac{y^n}{n!}\right)\frac{x^n}{n!} = e^y \sum_{n=0}^{\infty}(-1)^n L_n(x)\frac{y^{n+1}}{(n+1)!}.$$

Here is the proof of this identity. In order to compute the left hand side in (38) we consider the exponential generating function for the numbers $\frac{1}{n!}$

$$F(z) = \sum_{n=0}^{\infty}\frac{1}{n!}\frac{z^n}{n!} = \sum_{n=0}^{\infty}\frac{z^n}{(n!)^2} = I_0(2\sqrt{z})$$

where $I_0(x)$ is the modified Bessel function of zero order (see (33)). According to (7)

(39) $$\sum_{n=0}^{\infty}\left(e^y - 1 - \frac{y}{1!} - \ldots - \frac{y^n}{n!}\right)\frac{x^n}{n!} = e^y \int_0^y e^{-t} I_0(2\sqrt{xt})dt.$$

At the same time, by using the Cauchy rule for multiplication of power series we compute

(40) $$e^{-t}I_0(2\sqrt{xt}) = e^{-t}\sum_{n=0}^{\infty}\left\{\frac{x^n}{n!}\right\}\frac{t^n}{n!} = \left(\sum_{n=0}^{\infty}(-1)^n\frac{t^n}{n!}\right)\sum_{n=0}^{\infty}\left\{\frac{x^n}{n!}\right\}\frac{t^n}{n!}$$

$$= \sum_{n=0}^{\infty}\left\{\sum_{k=0}^{n}\binom{n}{k}\frac{(-1)^{n-k}x^k}{k!}\right\}\frac{t^n}{n!}.$$

The Laguerre polynomials $L_n(x)$ have the binomial representation

$$L_n(x) = \sum_{k=0}^{n} \binom{n}{k} \frac{(-1)^k x^k}{k!}$$

and from (39) and (40) we obtain the identity

$$\sum_{n=0}^{\infty} \left( e^y - 1 - \frac{y}{1!} - \cdots - \frac{y^n}{n!} \right) \frac{x^n}{n!} = e^y \int_0^y \left\{ \sum_{n=0}^{\infty} L_n(x) \frac{(-1)^n t^n}{n!} \right\} dt .$$

Now term by term integration leads to (38).